\documentclass[12pt]{article}
\usepackage{bm}        
\usepackage{amssymb,amsmath,amsopn,amsfonts,graphicx,xcolor}
\usepackage{algorithm,algpseudocode}
\newcommand{\beq}{\begin{equation}}
\newcommand{\eeq}{\end{equation}}
\newcommand{\beqn}{\begin{eqnarray}}
\newcommand{\eeqn}{\end{eqnarray}}
\newcommand{\R}{\mathbb R}

\begin{document}
	\title{Stochastic learning control of inhomogeneous quantum  ensembles}
	\author{Gabriel Turinici \\
		IUF - Institut Universitaire de France \\ 
		CEREMADE, Universit\'e Paris Dauphine - PSL Research University}
	
	\date{Oct 2019}
	\maketitle
	\begin{abstract}
In quantum control, the robustness with respect to uncertainties in the system's parameters or driving field characteristics is of paramount importance 
and has been studied theoretically, numerically and experimentally. 
We test in this paper  stochastic search procedures (Stochastic gradient descent and the Adam algorithm) that sample, at each iteration, from the distribution of the parameter uncertainty, as opposed to previous approaches that use a fixed grid. We show that both algorithms behave well with respect to benchmarks and discuss their relative merits. In addition the methodology allows to address high dimensional parameter uncertainty; we implement numerically, with good results, a 3D and a 6D case.


\end{abstract}
\maketitle

\section{Introduction}

Quantum control is a promising technology with many applications ranging from NMR~\cite{Glaser98}
to quantum computing~\cite{PhysRevLett.102.080501} and laser control of quantum dynamics~\cite{brifrabitz10}. The controlling field encounters many molecules which although identical in nature may interact differently with the incoming field because of e.g., different Larmor frequencies or rf attenuation factors (in NMR spin control or quantum computing, see~\cite{li-khaneja-06bis,SKINNER20038,PhysRevA.29.1419,bookquantumcomputation,PhysRevApplied.4.024012,kosut}), different spatial profile (see~\cite{rabitz:hal-00536535}) or other parameters (see~\cite{sugny_timeoptimal,praRabitz14ensemble,PhysRevA.100.022302}).
For obvious practical reasons, it is of paramount importance to ensure that the control quality is robust with respect to this heterogeneity.
 Thus the quantum control problem involves a unique set of driving fields $u(t) \in \R^L$, 
the same for all molecules in the ensemble, however each molecule is described by a set of parameters  $\theta \in \Theta \subset \R^d$ and the control outcome depends on both $u$ and $\theta$; the goal can be expressed as the maximization of the control quality averaged over $\theta$.
A different view is when the variability is not due to the presence of many different molecules but when uncertainties in the 
control implementation require to devise a field robust to fluctuations in those parameters.	

A first natural question is whether this is at all possible, i.e., if a single field can drive several distinct molecules to a common target; the answer is given by the theory of ensemble control controllability, see~\cite{0305-4470-37-1-019,Belhadj201523,li-khaneja-06bis,beauchard-bloch,borzibook} and is in general positive. However the theory does not explain how to find the control (except under specific regimes, see ~\cite{augier}). 
To do so, different algorithms have been proposed: the pseudo-spectral approach of Li et al. \cite{LiRuthsJCP09,LiJCP11,LiPNAS11} consider spectral and/or polynomial representations of the control problem in 2D ($d=2$); Wang considers iterative procedures based on sampling~\cite{WangAutomaticaEnsemble18}; the learning approach of Chen et al.~\cite{praRabitz14ensemble} 
and Kuang et al.~\cite{Kuang18} (the latter in the context of time-optimal control) consider a fixed uniform grid over the inhomogeneous parameter space and was tested for $d=2$. 
Finally, Wu and al.~\cite{robust3D} find robust controls using uniform grids in 2D and 3D ($d=3$).

In all these works there is always a fixed grid (or fixed sampling) involved when the control is searched. The rationale behind this idea 
is that a fixed grid makes the search more stable and a good choice of the grid is enough to describe efficiently the mean performance of the control over the parameter space in the spirit of a quadrature formula for the average over $\theta$. This is coherent with results from the approximation theory which inform that convergence is of order $e^{-\sqrt[d]{\mathcal N}}$, 
with respect to the number ${\mathcal N}$ of grid points; however the same formula indicates a bad scaling with respect to $d$. To address this 
{\it curse of dimensionality} and also explore the nature of the search landscape, 
we take here a different view: at each control iteration we use a new sampling in the spirit of Monte Carlo methods 
(see~\cite[Section 7.7]{nr3}) 
for computing high dimensional
 integrals. This will induce slight oscillations in the average but has the advantage to cover the space $\Theta$ of inhomogeneity even in high dimensions $d$. 
 	A similar approach has been tested independently in a very recent work by Wu and al.~\cite{Rebing_robust19} for a two-dimensional example and promising results were obtained; see Section~\ref{sec:algorithms} for comments on the differences between the two approaches.
The procedure we propose is detailed in the next section and the numerical results are the object of Section~\ref{sec:numerical}.

\section{Algorithms for ensemble quantum control}

We consider a control $u(t)=(u_1(t),...,u_L(t) ) \in \R^L$ acting on a molecule part of a larger ensemble. Each molecule is completely characterized by some inhomogeneity parameter $\theta \in \Theta \subset \R^d$ obeying a distribution law $P(\theta)$ on $\Theta$ (which can be the uniform distribution or any other). All molecules are subjected to the same control $u(t)$ during the time interval $[0,T]$ in order to reach some target.

\subsection{Evolution equations}

The dynamics of each molecule in the sample is governed by the Hamiltonian
$H(\theta,u) = H_0(\theta)+ \sum_{\ell=1}^{L} u_\ell(t) H_\ell(\theta)$ 
through the Schr\"odinger equation:
\beq
i \frac{d}{dt} \psi(t;\theta) = H(\theta,u) \psi(t;\theta),
\eeq
where $\psi$ is the wave-function of the molecule 
(here and below we set $\hbar=1$). 
Of course, $\psi$ depends on $u$ but for notational convenience we omit to write explicitly this dependence from now on. Once a finite dimensional basis $\{ |j\rangle, j=1,..., N\}$ is chosen, 
the state of the quantum system can be represented as 
\begin{equation} \label{eq:psi}
|\psi(t;\theta)\rangle= \sum_{j=1}^{N}c_{j}(t;\theta)|j\rangle. 
\end{equation}
Denoting
$C(t;\theta)=(c_{0}(t;\theta),...,c_{N}(t;\theta))^{T}$ the vector of coefficients $C$ satisfies the equation: 
\beq \label{eq:Ct}
\frac{d}{dt} C(t;\theta) = X(\theta,u) C(t;\theta),
\eeq
where $X$ is the representation of the Hamiltonian $H$ (including the $1/i$ factor) in the basis
$|j\rangle$, $j=1,..., N$.

Note that same setting also applies to non-linear Hamiltonians e.g. Bose-Einstein condensates (nonlinearity in $\psi$), 
or high order control terms~\cite{Dion99,Grigoriu09} (nonlinearity in $u$).

The quantum system can also be described in terms of a density matrix $\rho(t;\theta)$;  
this matrix is expressed in some basis for operators. 
Same happens when the molecule is coupled to a bath or when relaxation phenomena are at work, 
see~\cite{ALLARD19987}; in both cases the coefficients of this expansion follow an equation similar to \eqref{eq:Ct}.

\subsection{Optimization by stochastic gradient descent and Adam algorithms} \label{sec:algorithms}

 The control goal is encoded as the minimization, with respect to $u$, of an error, or "loss" functional  
 ${\mathcal L}(u,\theta)$ 
depending on the control $u$ and the Hamiltonian parameters $\theta$. 
When all the ensemble is considered, the following loss functional is to be minimized:

\beq \label{eq:J}
{\mathcal J}(u) = \int_\Theta {\mathcal{L}}(u,\theta) P(d \theta).
\eeq

The stochastic optimization algorithms described below construct an iterative process in order to find the $u$ that minimizes \eqref{eq:J}.

Historically the first to be considered, the stochastic gradient descent algorithm~\cite{SGD} (henceforth called SGD) consists in the following procedure:

\begin{algorithm}[H]
\caption{SGD} \label{algo:SGD}	
\begin{algorithmic}[1]

\State Choose a learning rate $\alpha > 0$, a mini-batch size $M>0$ and the initial control $u^0$.
		
\State Set iteration counter $k=0$.
		
\Repeat 
		
\State	\label{item:sgd1} 	Draw $M$ independent parameters $\theta_1^k, ..., \theta_M^k$ from the distribution $P(\theta)$ and compute the approximation 
		$g^k := \frac{1}{M}\sum_{m=1}^{M} \nabla_{u} \mathcal{L}(u^k;\theta_m^k) $ of the gradient 
		$\nabla_u{\mathcal J}(u^k)$ of ${\mathcal J}(\cdot)$
		at $u^k$.
		
\State 	set $u^{k+1} = u^k - \alpha g^k$ and $k=k+1$.

\Until	some stopping criterion is satisfied.
\end{algorithmic}
\end{algorithm}

In order to accelerate the convergence of the SGD algorithm, several improvements have been proposed (see~\cite{ruder2016overview}) among which the Adam~\cite{kingma2014adam} variant which proved to be one of the most efficient and very scalable. The difference between Adam and SGD is that Adam uses a different learning rate for each parameter which is tuned as follows: when the uncertainty in the gradient is large the learning rate is taken to be small and contrary otherwise. In order to have a robust estimation for the gradient (in absolute value) a Exponential Moving Average is computed on the fly (see below). 
It can be described as:
\begin{algorithm}[H]
	\caption{Adam} \label{algo:Adam}	
	\begin{algorithmic}[1]
\State Choose the learning rate $\alpha > 0$, the EMA parameters $\beta_1$ and $\beta_2$, the mini-batch size $M>0$,
the epsilon $\varepsilon>0$
and the initial control $u^0$.
\State Set iteration counter $k=0$, first moment estimate $\mu=0$, second moment estimate $v=0$.
\State \label{item:adam1} Set $k=k+1$.
\Repeat 
\State	Draw $M$ independent parameters $\theta_1^k, ..., \theta_M^k$ from the distribution $P(\theta)$ and compute the approximation 
$g^k := \frac{1}{M}\sum_{m=1}^{M} \nabla_{u} \mathcal{L}(u^{k-1};\theta_m^k) $ of the gradient 
$\nabla_u{\mathcal J}(u^{k-1})$ of ${\mathcal J}(\cdot)$
at $u^{k-1}$.
\State Compute the moving averages $\mu^k := \beta_1 \mu^{k-1} + (1-\beta_1) g^{k}$,
$v^k := \beta_2 v^{k-1} + (1-\beta_2) |g^{k}|^2$.
\State   \label{item:adambiasstep} Compute bias-corrected moment estimates: 
$\hat{\mu}^k = \mu^k / (1-(\beta_1)^k)$, $\hat{v}^k = v^k / (1-(\beta_2)^k)$.
\State set $u^{k} = u^{k-1} - \alpha \hat{\mu}^k / (\sqrt{\hat{v}^k}+\varepsilon)$.
\Until some stopping criterion is satisfied.
	\end{algorithmic}
\end{algorithm}

	The momentum algorithm used in~\cite{Rebing_robust19} can be seen as being halfway  between SGD and Adam; it is formally a special case of the Adam algorithm for $\beta_1 = \lambda$, $\beta_2=1$, $v^0=1$ and no bias correction step~\ref{item:adambiasstep} (that is 
	$\hat{\mu}^k = \mu^k$,  $\hat{v}^k = v^k $). In practice the numerical results are very similar and point in the same direction; in particular we expect that the momentum algorithm is also relevant to high dimensional robust control problems.

\section{Numerical results} \label{sec:numerical}

We test the performance of the algorithms in Section~\ref{sec:algorithms} for several benchmarks from the literature (or that generalize cases from the literature).

	 In sections~\ref{sec:Rabitz2level} and~\ref{sec:Rabitz2level} we compare the SGD algorithm with a fixed grid
	sampling method from the literature. Then in section   sections~\ref{sec:Li} and~\ref{sec:Wang}  we compare the SGD wih the Adam 
	algorithm and in Section~\ref{sec:histograms} we draw further conclusions concerning stochastic optimization.

In the situations considered below, the goal is to maximize the so-called {\it fidelity} denoted $\mathcal{F}(u;\theta)$. 
For sections~\ref{sec:Rabitz2level} and 
  \ref{sec:Rabitz3level} this has the formula 
  $\mathcal{F}(u;\theta)= |\langle C(T;\theta), C_{\text{target}} \rangle |$ where $C_{\text{target}}$ is a prescribed target state. But this expression is not differentiable everywhere and
  numerically it is easier to replace it with its square.
   Moreover, to express the problem as a minimization, a $-1$ multiplicative constant is introduced and $1$ added to the result in order to have it positive.
   So the cost functional 
  ${\mathcal J}$ will be the mean, over $\theta \in \Theta$ of the error in the fidelity squared as in formula~\eqref{eq:JRabitz2D}.
  On the contrary, when the fidelity is more well behaved  as in section~\ref{sec:Li} where $\mathcal{F}(u;\theta)=c_4(T,\theta)$
  or in section~\ref{sec:Wang}  where $\mathcal{F}(u;\theta)=c_6(T,\theta)$ the square operation is useless and the cost functional has the form in 
  \eqref{eq:J3D} or \eqref{eq:J6D}. However, in  all sections, we will plot the
  error in the fidelity itself; the reason why not plotting the fidelity (instead of the error) is that the error can be very small 
    (as in Section~\ref{sec:Rabitz2level}) and the results are more visible on a logarithmic scale.
 Note that in some cases the best control cannot attain the target with 100\% quality (even for a single molecule). However, 
	for any given value of the parameter $\theta$, the best attainable performance is 
	known (see~\cite{Glaser98,Turinici20011,KhanejaPnas03}) and is denoted $F_{max}(\theta)$. 
	We will therefore consider the fidelity relative to $F_{max}(\theta)$.
 In all cases the error is computed as the average over $M_{test}=300$ random independent parameters 
$\theta_1^{test}$, $\theta_2^{test}$, ..., $\theta_{M_{test}}^{test}$ drawn (once for all) from the distribution $P(\theta)$ and has
 the following expression:
\begin{equation} \label{eq:deferrorplot}
\frac{1}{M_{test}}\sum_{k=1}^{M_{test}} \left( 1-\frac{ \mathcal{F}(u;\theta_k^{test})   }{F_{max}(\theta_k^{test})} \right).
\end{equation}
For sections~\ref{sec:Rabitz2level} and 
\ref{sec:Rabitz3level} we will also plot the max relative error: 
\begin{equation} \label{eq:defmaxerrorplot}
\max_{k=1, ..., M_{test}} \left( 1-\frac{ \mathcal{F}(u;\theta_k^{test})   }{F_{max}(\theta_k^{test})} \right).
\end{equation}

Finally, in order to compare our algorithm with those from the literature, we take as indicator of the numerical effort the number of gradient $\nabla_{u} \mathcal{L}(u;\theta)$ evaluations; for instance one iteration of SGD or Adam algorithms count as $M$ gradient evaluations. 
In all situations we used for the Adam algorithm the standard values 
$\beta_1=0.9$, $\beta_2=0.999$, $\varepsilon=10^{-8}$.

\subsection{Two level inhomogeneous ensemble}  \label{sec:Rabitz2level}

Consider an ensemble of spins as in~\cite[section III.]{praRabitz14ensemble}. 
The spins have different Larmor frequencies $\omega$ in the range $[0.8,1.2]$ and the controls ($L=2$) have multiplicative inhomogeneity $\epsilon \in [0.8,1.2]$; 
we set $\theta= (\omega,\epsilon)$ and
with the previous notations the dynamics corresponds to the equation:
\beq
{\scriptsize
\begin{pmatrix}
	\dot{c}_{1}(t;\theta) \\
	\dot{c}_{2}(t;\theta) \\
\end{pmatrix}=
 \begin{pmatrix}
	0.5\omega i & 0.5\epsilon(u_2(t)- i u_1(t))  \\
	-0.5\epsilon(u_2(t)- i u_1(t))  & -0.5 \omega i \\
\end{pmatrix}
\begin{pmatrix}
	c_{1}(t;\theta) \\
	c_{2}(t;\theta) \\
\end{pmatrix},
}\eeq
	where $c_1$, $c_2$ are the coefficients of the wavefunction of the spin system in the canonical basis, as detailed in 
equation~\eqref{eq:psi}.

The initial state of each member of the quantum ensemble is set to
$|\psi_{0}\rangle=|0\rangle$; i.e., $C_{0}=(1,0)^T$, 
and the goal is to reach the target
state $|\psi_{\text{target}}\rangle=|1\rangle$; i.e.,
$C_{\textrm{target}}=(0,1)^T$. 
The objective is encoded as the requirement to minimize:
\beq \label{eq:JRabitz2D}
J(u) = \frac{1}{2} \left(1-\int_\Theta |\langle C(T;\theta), C_{\text{target}} \rangle |^2 P(d\theta) \right).
\eeq

Here $F_{max}(\theta)=1$.
The total time is $T=2$ is divided into  $Q=200$ time steps, of length $\Delta t=T/Q=0.01$ each.
The initial choice for the control $u$ is
$u^{k=0}(t)=\{u^{0}_{1}(t)=\sin t, u^{0}_{2}(t)=\sin t\}$.

\begin{figure*}
	\includegraphics[width=\textwidth]{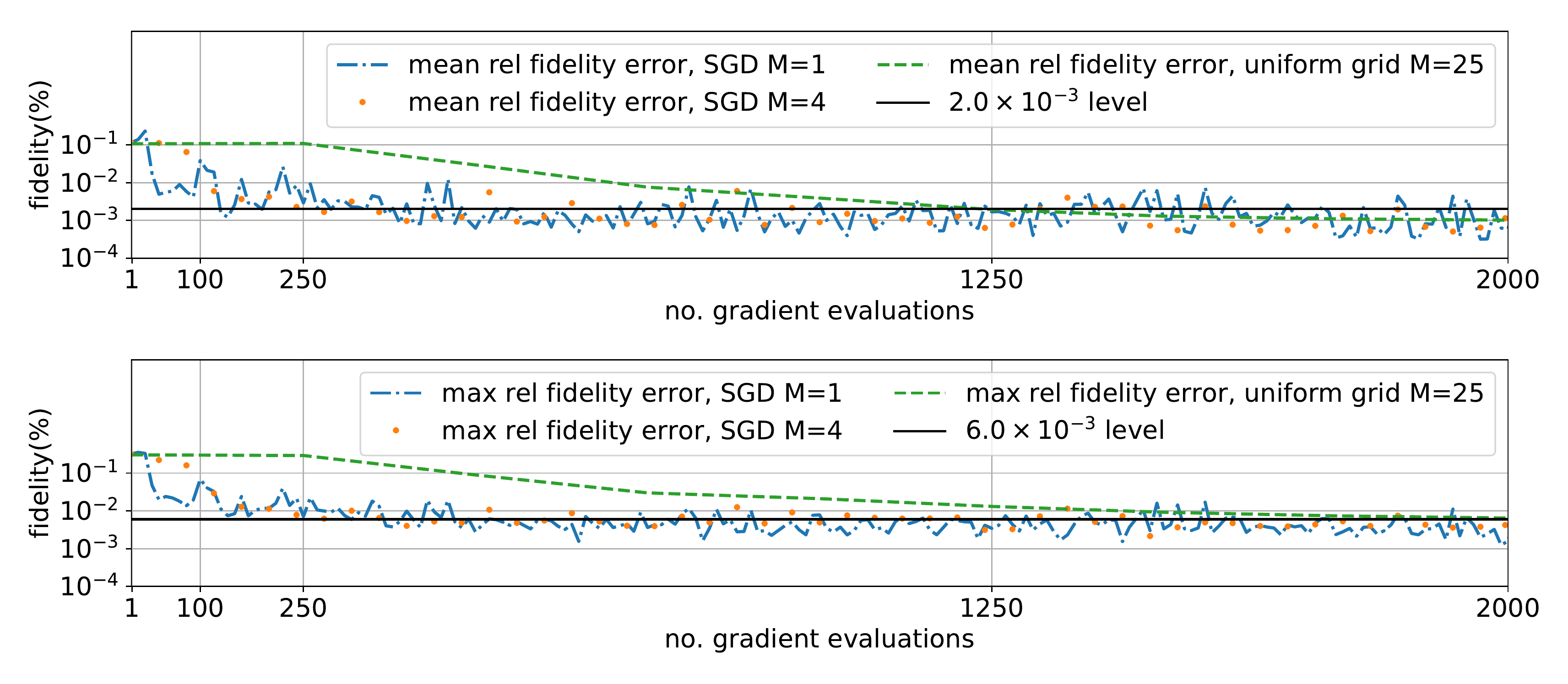}%
	\caption{\label{fig:Rabitz2d} 
		Convergence for the numerical case in Section~\ref{sec:Rabitz2level}. 
		Top image:  mean fidelity error	(as defined in equation~\eqref{eq:deferrorplot}). 
		Bottom image: maximum (over the sample) fidelity error (as defined in equation~\eqref{eq:defmaxerrorplot}).
		We consider three simulations: 
		a fixed uniform 2D grid ($M=25$) as in~\cite[section III.A.]{praRabitz14ensemble} and the SGD algorithm with $M=1$ and $M=4$. 
		This SGD   
		converges about $5$ times faster: the mean fidelity error of $2.0\times 10^{-3}$ is obtained after $1250$ gradient evaluations of the fixed grid 
		algorithm and
		after  $250$ evaluations of the SGD algorithm with $M=1, 4$. Same for other levels of errors.
	}
\end{figure*}

Several mini-batch sizes $M=1, 4, 8, 16$ and $32$ are tested and compared 
with implementation in~\cite[section III.A.]{praRabitz14ensemble} where a 2D uniform grid of $5 \times 5$   values for $\theta$  is chosen.
In all cases very good convergence results are attained. We plot in Figure~\ref{fig:Rabitz2d} the results for $M=1$, $M=4$  
relative to the convergence with the uniform $5 \times 5$  grid. 
In all cases ($M=1, 4$, uniform grid) we set $\alpha=500$; note that the learning rate  
$\alpha$ was optimized to obtain the best possible results for the 
fixed grid algorithm and indeed the results are better than those in~\cite[section III.A.]{praRabitz14ensemble}.
But similar conclusions are reached for any value of $\alpha$. An acceleration by a factor of $5$ is obtained
 for both $M=1$ and $M=4$, essentially due to the fact that each SGD iteration uses only 
$M$ gradient evaluations.
 Note that the SGD algorithm oscillates
but these oscillations can be cured by lowering $\alpha$ (or stopping the search) as soon as a good result is obtained.
The question of which is the best choice among $M=1$ and $M=4$ is a matter of striking a balance between speed and uncertainty: for
$M=4$ the convergence is slightly slower but oscillations are diminished. This behavior is observed, to a larger or lesser extent, in all test cases.

Note that in order to compare our learning rate $\alpha$ (for the fixed uniform grid) with that in \cite[section III.A.]{praRabitz14ensemble} a
	multiplicative factor of $\Delta t /2$ has to be introduced because our gradient (see Appendix~\ref{sec:gradientcomput}) contains an extra $\Delta t$ factor
and the coefficient $1/2$. Thus one should transform $\alpha=500$ to $1/2*0.01*500 = 2.5$ to compare with $0.2$ used in~\cite{praRabitz14ensemble}.

\subsection{A three level $\Lambda$ atomic ensemble} \label{sec:Rabitz3level}

In this section we test a $\Lambda$ atomic ensemble from~\cite[Section IV]{praRabitz14ensemble}
which can be written as a  $3$-level system with the following dynamics:

\begin{eqnarray}\label{3level-element1}
\!\!\!\!\!\!\!\!  
{\small
\begin{pmatrix}
\dot{c}_{1}(t;\theta) \\
\dot{c}_{2}(t;\theta) \\
\dot{c}_{3}(t;\theta) \\
\end{pmatrix}=
\begin{pmatrix}
-1.5\omega i & 0 & -i\epsilon u_{2}(t)  \\
0 & -\omega i  & -i\epsilon u_{1}(t) \\
-i\epsilon u_{2}(t) & -i\epsilon u_{1}(t) & 0 \\
\end{pmatrix}
\begin{pmatrix}
c_{1}(t;\theta) \\
c_{2}(t;\theta) \\
c_{3}(t;\theta) \\
\end{pmatrix} \!\!,
}\end{eqnarray}
where $\omega$ and $\epsilon$ 
have uniform distributions in  $[0.8, 1.2]$
	and $c_1$, $c_2$, $c_3$ are the coefficients of the wavefunction of the spin system in the canonical basis, as detailed in 
	equation~\eqref{eq:psi}.

The objective
is to find a control  $u(t)=(u_1(t),u_2(t))$ which drives
all the inhomogeneous members from
$|\psi_{0}\rangle=\frac{1}{\sqrt{3}}(|1\rangle+|2\rangle+|3\rangle)$
(i.e.,
$C_{0}=(\frac{1}{\sqrt{3}},\frac{1}{\sqrt{3}},\frac{1}{\sqrt{3}})$)
to $|\psi_{\text{target}}\rangle=|3\rangle$ (i.e.,
$C_{\text{target}}=(0,0,1)$); the objective is encoded as the minimization of~\eqref{eq:JRabitz2D}.
Here $F_{max}(\theta)=1$.

\begin{figure*}
	\includegraphics[width=\textwidth]{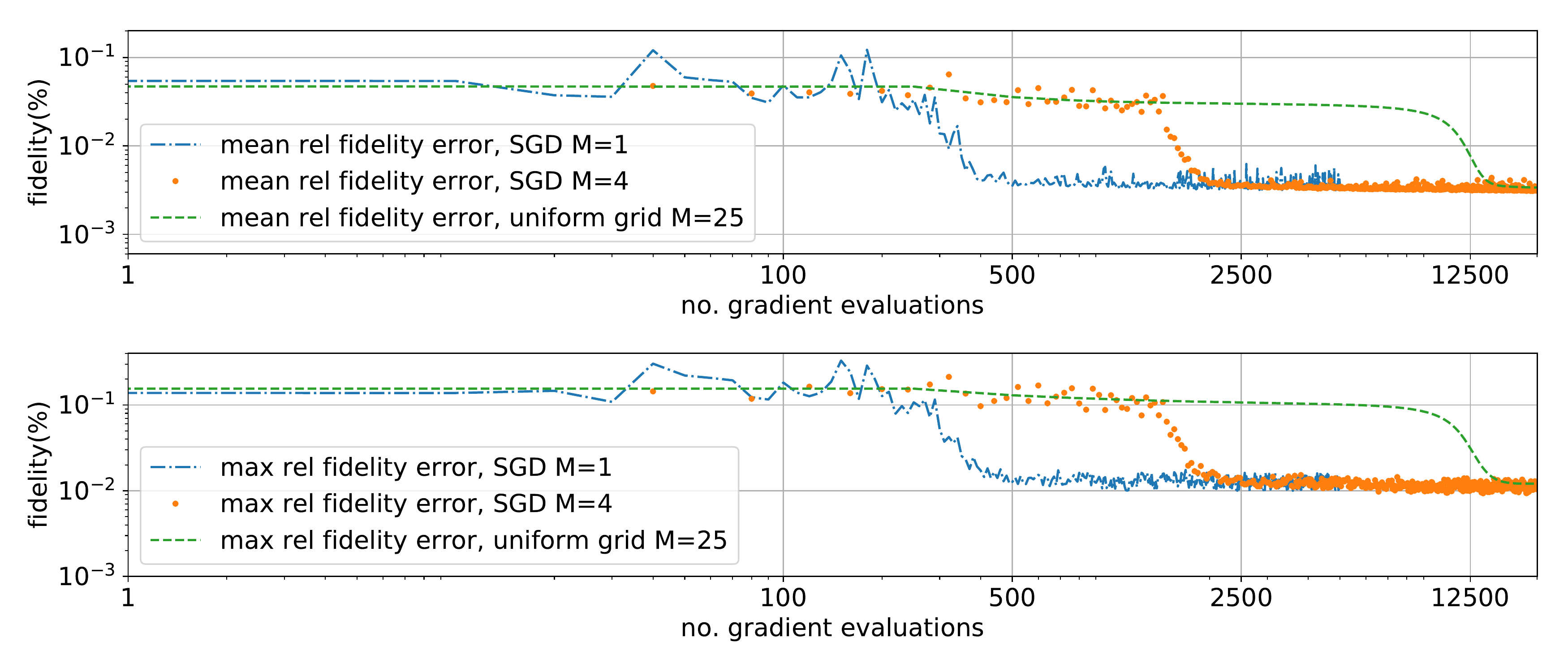}%
	\caption{\label{fig:Rabitz3d} 		
		Convergence for the numerical case in Section~\ref{sec:Rabitz3level}. 
		Top image:  mean fidelity error as defined in equation~\eqref{eq:deferrorplot}. 
		Bottom image: maximum (over the sample) fidelity error (as defined in equation~\eqref{eq:defmaxerrorplot}).
		We consider two algorithms: 
		a fixed uniform 2D grid ($M=25$) as in~\cite[section IV.]{praRabitz14ensemble} and the SGD algorithm with 
		$M=1$ and $M=4$. 
		This latter approach 
		converges about $7$ times faster: the convergence settles in after $17'500$ 
		gradient evaluations of the fixed grid algorithm compared with cca. $2'500$ evaluations of the SGD algorithm. 
		This acceleration factor is even more important for $M=1$  but at the price of larger oscillations.
	}
\end{figure*}

 We plot in Figure~\ref{fig:Rabitz3d} the results for $M=1$ and $M=4$  
 relative to the convergence with an uniform grid as in~\cite[section IV.]{praRabitz14ensemble}.
In all cases ($M=1, 4$, uniform grid) we set $\alpha=100$.
The acceleration factor is around $7$ for $M=4$ and even larger for $M=1$ (but at the price of larger oscillations too).

\subsection{A 3D example: two spin systems without cross-correlated relaxation} \label{sec:Li}

As argued before, methods from the literature may have difficulties to address high dimensional parameters, 
and often limit to two dimensional ($d=2$) inhomogeneity 
(see however~\cite[Sec. V.B]{robust3D} for a 3D case).
In order to test the full power of our method, we consider two situations that extend cases treated in the literature
but have never been treated before. 
The first test is a three dimensional ($d=3$) example
which addresses the coherence transfer between two spins without cross-correlated relaxation, taken from~\cite[Section III.B.1. eq(15)]{LiJCP11} (but with an additional inhomogeneity dimension). 
An example of such a system is an isolated hetero-nuclear spin system composed of two coupled spins $1/2$ 
corresponding to atoms $\ ^1\!H$ and $\ ^{15}\!N$.
For a general treatment of the relaxation terms and the formulation of this equation see~\cite{ALLARD19987}.
The spins
display control inhomogeneity described by the parameter $\epsilon$ as above but there is also variation
in the relaxation rate and coupling constant, which, denoting $\theta=(\epsilon,J,\xi)$ results in the dynamical system:
\begin{eqnarray}\label{Li2spinnocrosscor}
\!\!\!\!\!\!\!\!  
{\small
\begin{pmatrix}
\dot{c}_{1}(t;\theta) \\
\dot{c}_{2}(t;\theta) \\
\dot{c}_{3}(t;\theta) \\
\dot{c}_{4}(t;\theta) \\
\end{pmatrix}=
\begin{pmatrix}
0 & - \epsilon u_1(t) & 0 & 0 \\ 
\epsilon u_1(t) & -\xi & -J & 0 \\ 
0 & J & -\xi & -\epsilon u_2 \\
0 & 0 & \epsilon u_2 & 0
\end{pmatrix}
\begin{pmatrix}
c_{1}(t;\theta) \\
c_{2}(t;\theta) \\
c_{3}(t;\theta) \\
c_{4}(t;\theta) \\
\end{pmatrix} \!\!.
}\end{eqnarray}
Let us denote by $I_{1x} = \sigma_x/2, I_{1y}= \sigma_y/2, I_{1z}= \sigma_z/2$ 
(here $\sigma_x$, $\sigma_y$ $\sigma_z$ are the Pauli matrices)
the spin operators 
corresponding to the first spin and 
 $I_{2x}, I_{2y}, I_{2z}$ the corresponding objects for the second spin. With the usual notations for the Kronecker products, 
$c_1 = \langle I_{1z} \rangle $, 
$c_2 = \langle I_{1x} \rangle $, 
$c_3 = \langle 2 I_{1y} I_{2z}\rangle $, 
$c_4 = \langle 2 I_{1z} I_{2z}\rangle $; 
	 the exact derivation
	 of this equation is beyond the scope of this work,  see~\cite{ALLARD19987,Glaser98,KhanejaPnas03} for details.	 
On the other hand also note that the dynamics is not reversible (relaxation is present) and the equations do not correspond to a unitary evolution.

The inhomogeneity  $\theta=(\epsilon,J,\xi)$ is uniformly distributed in 
$\Theta= [0.9,1.1]\times [0.5,1.5]\times [0,2]$. The final time $T=7\pi/6$ is discretized with 
$Q=200$ uniform time steps. The control is initialized as  before. The initial state is 
encoded as $c_{0}=(1,0,0,0)$
and the target is to minimize the 
three-dimensional integral:
\begin{equation} \label{eq:J3D}
{\mathcal J}(u)= 1- \int_\Theta c_4(T;\theta) P(d\theta).
\end{equation}

\begin{figure}
	\includegraphics[width=3.5in]{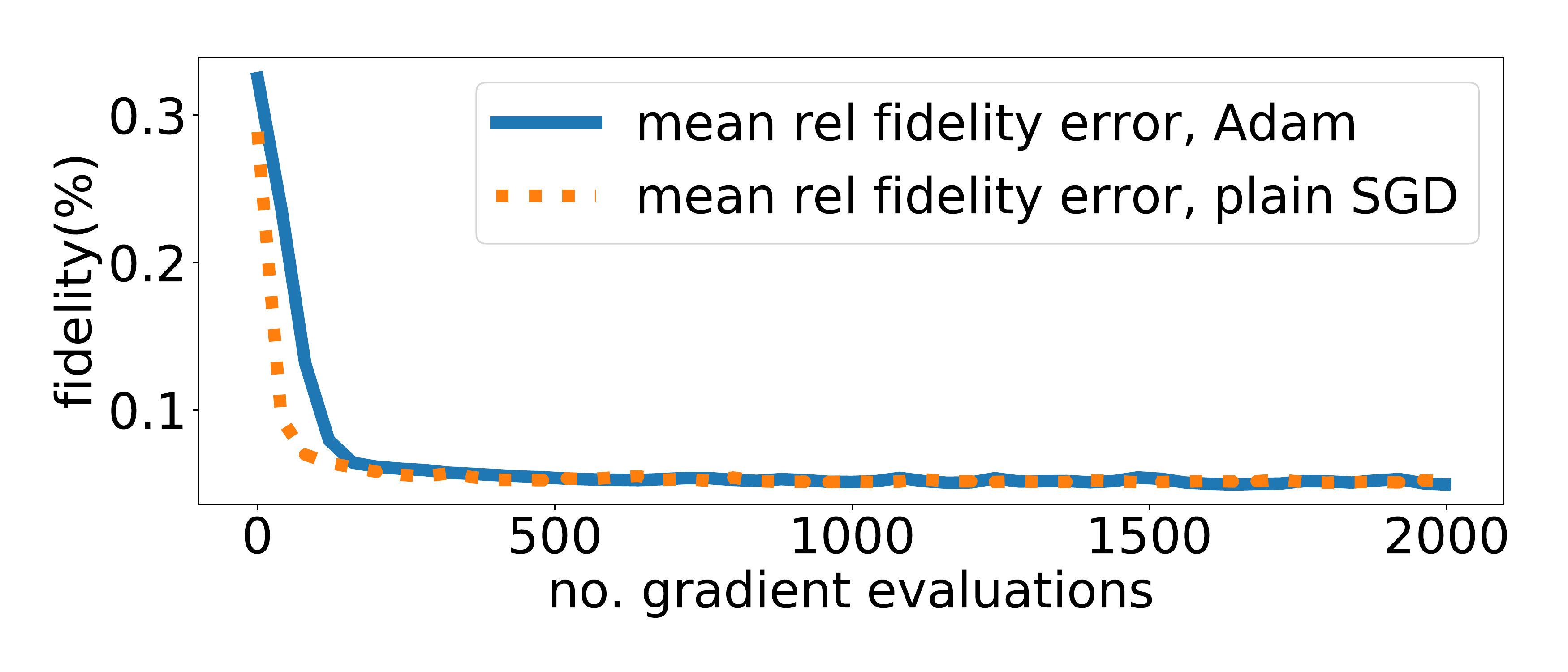}%
	\caption{\label{fig:Li3p1} 
		Convergence for the numerical case in Section~\ref{sec:Li}.
		The quantity plotted is given in equation~\eqref{eq:deferrorplot}.		 
		We set $M=4$; for the SGD algorithm we choose $\alpha=10.0$ and for the Adam algorithm we set $\alpha=0.01$.
		The continuous ($-$) and dotted ($\cdot$)  curves stand for the mean fidelity errors
		of the plain SGD and Adam algorithm respectively; the convergence is similar and a 
		$95\%$ mean target relative fidelity 
		(or equivalently $5\%$ mean target relative fidelity error)
		is obtained.
		For the controls see Figure~\ref{fig:Li3p2}.
	}
\end{figure}

\begin{figure}
	\includegraphics[width=3.5in]{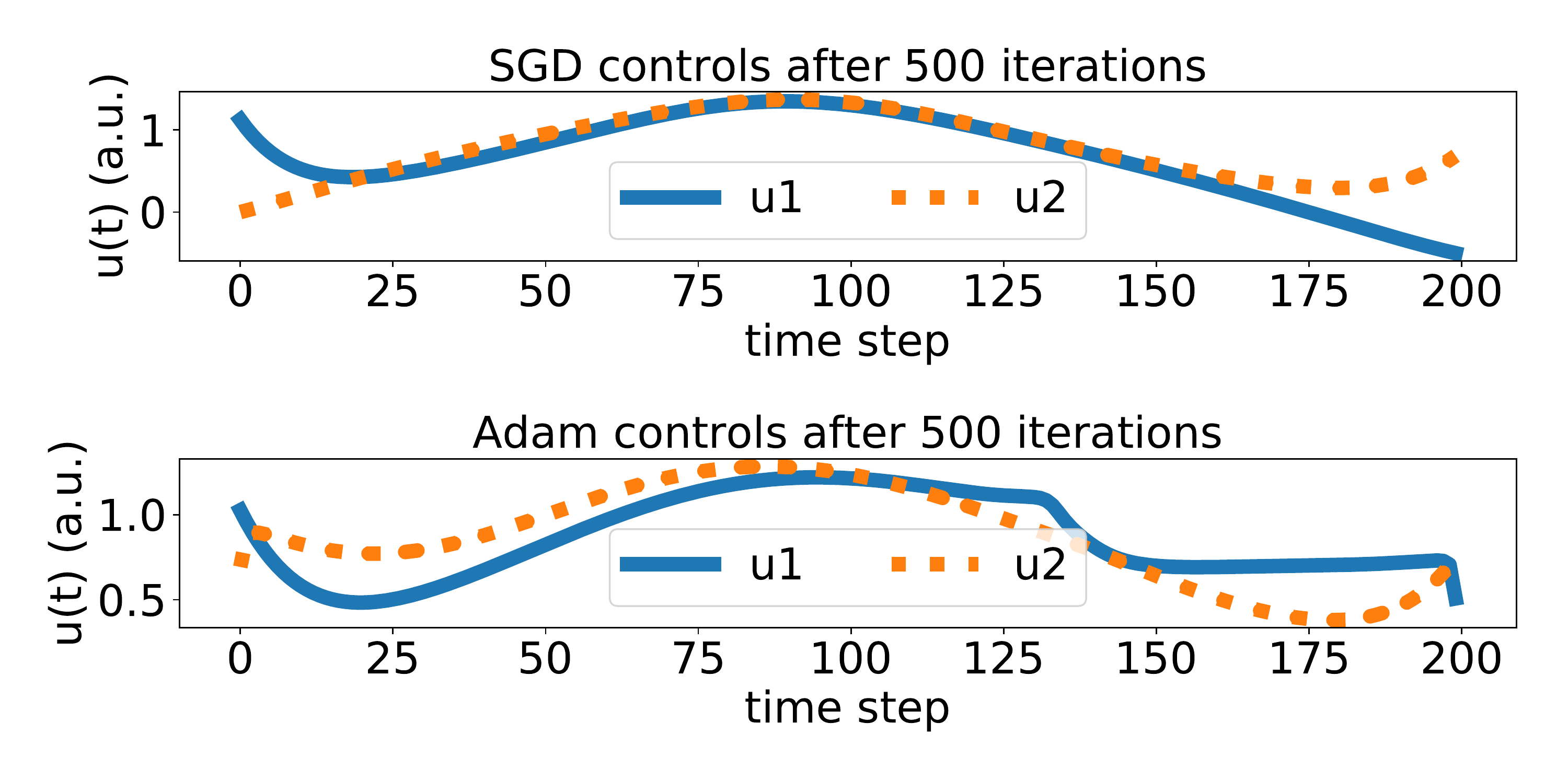}%
	\caption{\label{fig:Li3p2} Converged controls for the SGD (up) and Adam (bottom) for the situation in 
		in Section~\ref{sec:Li} (for the convergence see Figure~\ref{fig:Li3p1}).  
		Controls obtained with the SGD algorithm are smoother than those from the Adam algorithm.
	}
\end{figure}

Recall that here the fidelity is $\mathcal{F}(u;\theta)=c_4(T,\theta)$; in this case 
(see~\cite{Glaser98,KhanejaPnas03}) $F_{max}(\theta)=\sqrt{1+(\xi/J)^2} -\xi/J$ (the worse performance being $-F_{max}(\theta)$).
The results are in Figures~\ref{fig:Li3p1} and \ref{fig:Li3p2}.
 Note that although for each $\theta$ 
 taken individually
 the figure $F_{max}(\theta)$ can be attained with a pair (recall $L=2$) of suitable control fields, 
 it is unknown whether a unique control pair exists ensuring $100\%$ (relative to $F_{max}(\theta)$) target yield simultaneously for all $\theta \in \Theta$. 
 In practice we did not find any, irrespective of the algorithm hyper-parameters such as $\alpha$, the  maximum number of iterations etc.; 
 we conclude on one hand that this ensemble is not $100\%$ simultaneously controllable and on the other hand that  
our procedure  improves significantly 
the robustness of the control with respect to  $\theta \in \Theta$ from an initial value of 
$67\%$ up to $95\%$. Note that the results from the literature (which for this case
 only consider 2 dimensional inhomogeneity) do not obtain $100\%$ control either (exact figure is not reported).
  
\subsection{A 6D example: two spin systems with cross-correlated relaxation} \label{sec:Wang}

We continue here to address new systems that previous methods could not treat. We consider an ensemble of two spin systems 
with cross-correlated relaxation as in
\cite[Section III.A.2.]{LiRuthsJCP09},  \cite[Section III.B.2 eq. (16)]{LiJCP11} and also
\cite[Example 3]{WangAutomaticaEnsemble18}, \cite{ALLARD19987}.

The spins 
display control inhomogeneity described by the parameters $\epsilon_1$ and $\epsilon_2$
and there is also variation
in the auto-correlated relaxation rate $\xi_a$,
the quotient $\xi_c/\xi_a$ of the cross-correlation relaxation rate $\xi_c$ with respect to the 
auto-correlated relaxation rate $\xi_a$ and finally, a dispersion in the Larmor frequencies of each spin.
Denoting $\theta=(\epsilon_1,\epsilon_2,\omega_1,\omega_2,\xi_a,\xi_c/\xi_a)$ $\in$ 
$\Theta=[0.9,1.1]^2 \times [0,1]^2 \times [0.75,1.25] \times [0.7,0.9]$, the dynamical system can be written:

\begin{eqnarray}\label{Li2spinwcrosscor}
\begin{pmatrix}
\dot{c}_{1}(t;\theta) \\
\dot{c}_{2}(t;\theta) \\
\dot{c}_{3}(t;\theta) \\
\dot{c}_{4}(t;\theta) \\
\dot{c}_{5}(t;\theta) \\
\dot{c}_{6}(t;\theta) \\
\end{pmatrix}=
\begin{pmatrix}
0 & - \epsilon_1 u_1(t) & \epsilon_2 u_2(t) &  0 & 0 & 0\\ 
\epsilon_1 u_1(t) & -\xi_a & \omega_1 & -J & -\xi_c & 0 \\ 
 -\epsilon_2 u_2(t) & -\omega_1 & -\xi_a & -\xi_c & J& 0 \\
0 & J & -\xi_c & -\xi_a &  \omega_2 & -\epsilon_2 u_2(t) \\
0 & -\xi_c & -J & -\omega_2 & -\xi_a & \epsilon_1 u_1(t)  \\
0 & 0 & 0& \epsilon_2 u_2(t) & -\epsilon_1 u_1(t) & 0 \\
\end{pmatrix}
\begin{pmatrix}
c_{1}(t;\theta) \\
c_{2}(t;\theta) \\
c_{3}(t;\theta) \\
c_{4}(t;\theta) \\
c_{5}(t;\theta) \\
c_{6}(t;\theta) \\
\end{pmatrix}.
\end{eqnarray}
The vector $C=(c_1,...,c_6)$ has real entries
 and, 
 with the same notations as in equation~\eqref{Li2spinnocrosscor},
  $c_1 = \langle I_{1z} \rangle $, 
 $c_2 = \langle I_{1x} \rangle $, 
 $c_3 = \langle I_{1y} \rangle $, 
 $c_4 = \langle 2 I_{1y} I_{2z}\rangle $, 
 $c_5 = \langle 2 I_{1x} I_{2z}\rangle $, 
 $c_6 = \langle 2 I_{1z} I_{2z}\rangle $.
 The relations are
 similar to that in section~\ref{sec:Li},
with the exception that there are two new entries $c_3$ and $c_5$  due to the presence of cross-correlation,
 see~\cite{Glaser98,KhanejaPnas03,ALLARD19987} for details of the derivation of the model; the dynamics is not reversible (relaxation is present) nor unitary.

We set $J=1$;  the total time $T=5$ is discretized with 
$Q=200$ uniform time steps. The control is initialized as before. The initial state is 
encoded as $c_{0}=(1,0,0,0,0,0)$
and the target is to minimize the six-dimensional integral:
\begin{equation} \label{eq:J6D}
{\mathcal J}(u)= 1- \int_\Theta c_6(T;\theta) P(d\theta).
\end{equation}

Recall that here the fidelity is $\mathcal{F}(u;\theta)=c_6(T,\theta)$.
In this case too, the best attainable performance for a single molecule is known (see~\cite{Glaser98,KhanejaPnas03}) and defined by 
$F_{max}(\theta)=\sqrt{1+\eta^2} - \eta$ where $\eta=\sqrt{\frac{\xi_a^2 - \xi_c^2 }{J^2 + \xi_c^2}}$.

\begin{figure}
	\includegraphics[width=3.5in]{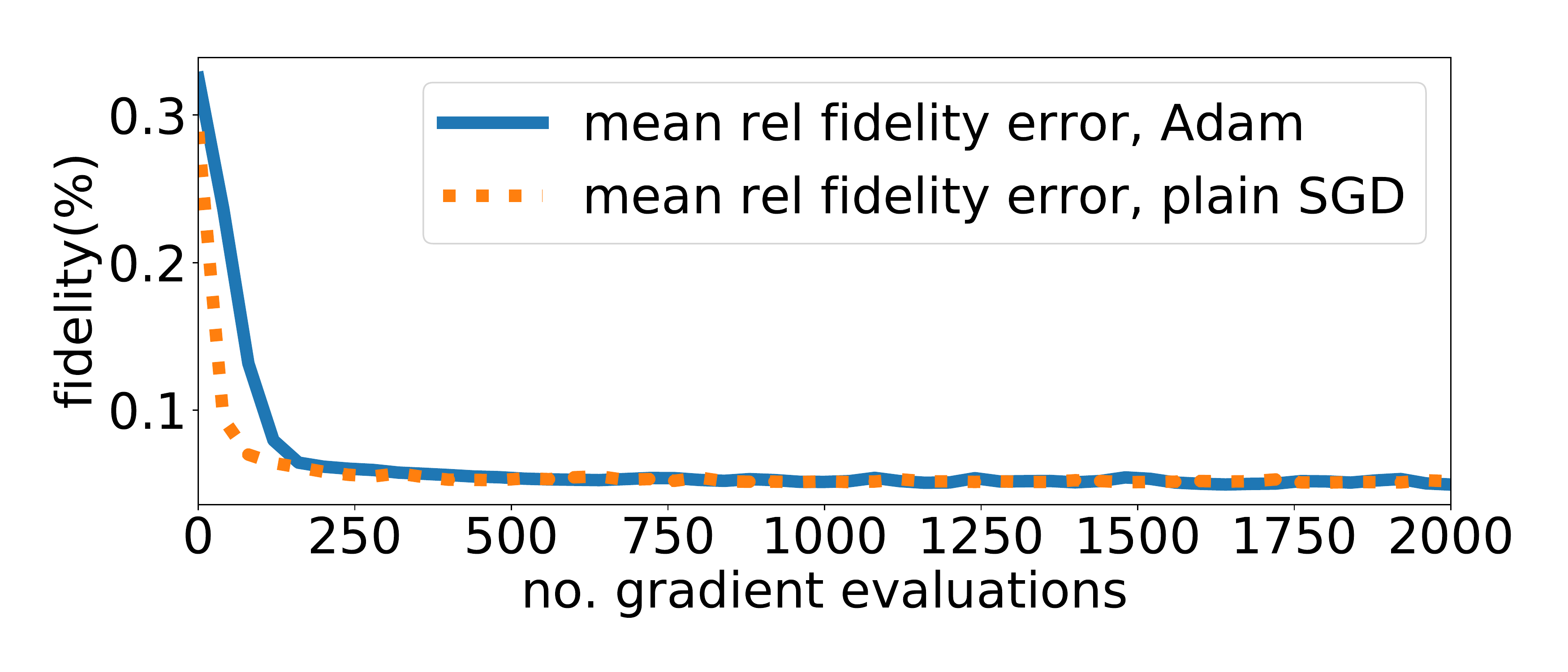}%
	\caption{\label{fig:Wang6p1} 
		Convergence for the numerical case in Section~\ref{sec:Wang}.
		The quantity plotted is defined as in the Figure~\ref{fig:Li3p1}.		 
		We set $M=4$; for the SGD algorithm we choose $\alpha=10.0$ and for the Adam algorithm we set $\alpha=0.01$.
		The continuous ($-$) and dotted ($\cdot$)  curves stand for the mean fidelity errors 
		of the plain SGD and Adam algorithm respectively; the convergence is similar and $91\%$ mean relative fidelity is obtained.
		For the controls see Figure~\ref{fig:Wang6p2}.
	}
\end{figure}

\begin{figure}
	\includegraphics[width=3.5in]{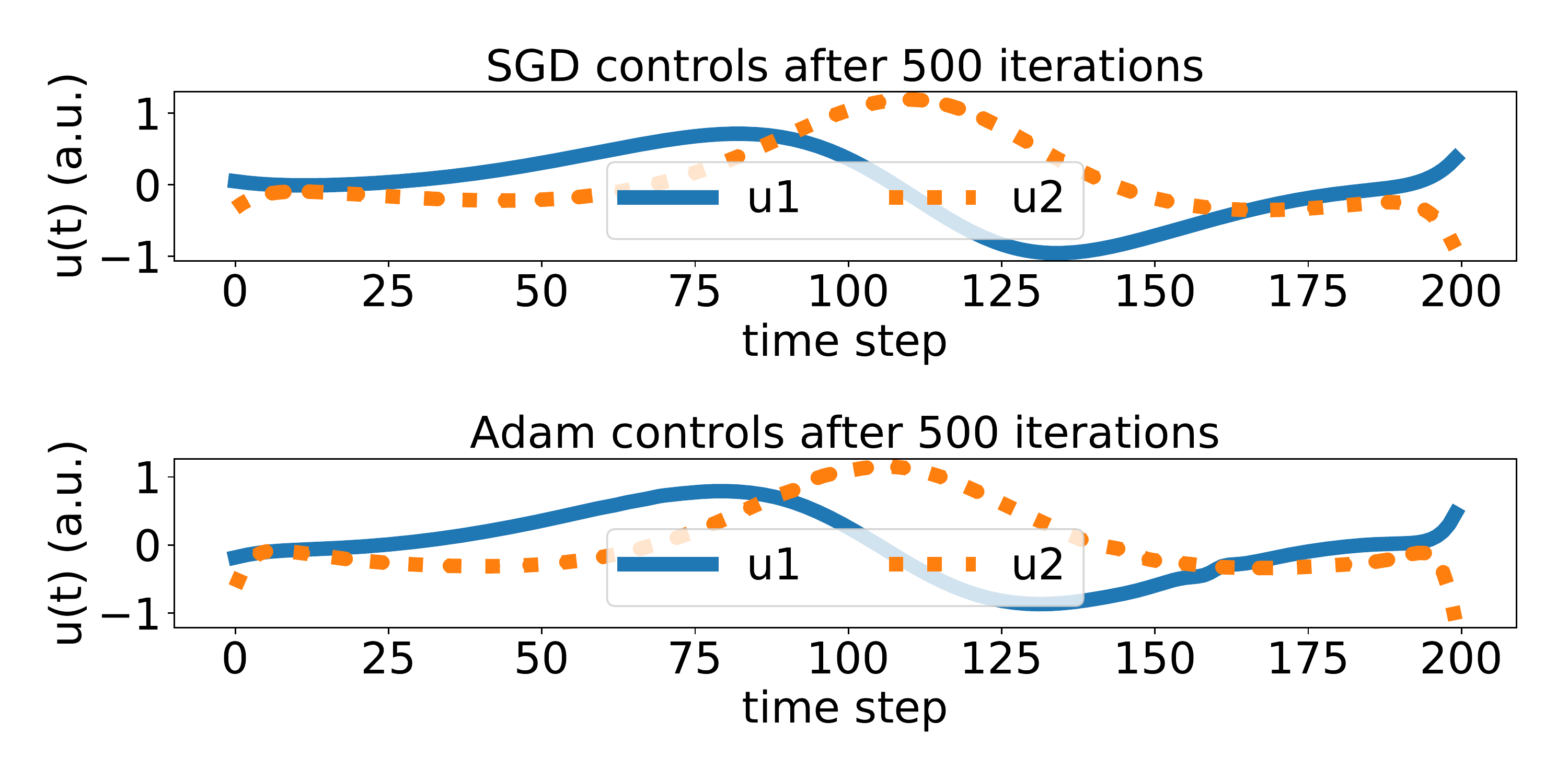}%
	\caption{\label{fig:Wang6p2} Converged controls for the SGD (up) and Adam (bottom) for the situation in 
		in Section~\ref{sec:Wang} (for the convergence see Figure~\ref{fig:Wang6p1}).  
		Controls obtained with the SGD algorithm are smoother than those from the Adam algorithm.
	}
\end{figure}

The simulation results are in Figures~\ref{fig:Wang6p1} and \ref{fig:Wang6p2}. Same conventions are kept as in the previous section (fidelity is relative to 
maximum attainable figure) and same considerations still apply: $100\%$ simultaneous controllability does not seem attainable but
significant improvement in the robustness is obtained ($91\%$ up from $-8\%$).

\subsection{Stochastic convergence behaviors} \label{sec:histograms}

The convergence of the stochastic algorithms can have two important regimes: 
\begin{enumerate}
	\item first, when all members of the ensemble can be simultaneously optimized to $100\%$; in our situation this is equivalent to simultaneous controllability. In this case convergence is "easier" because  it is "enough" to follow the gradient for each parameter value in order to converge; at convergence all gradients (as distribution with respect 
	to $\omega$), will collapse to (in practice will be close to) a Dirac mass.
	\item secondly, when members of the ensemble cannot be simultaneously optimized; in this case, reaching full control for some $\theta$ value will harm the quality of some other parameter values $\theta' \neq \theta$. At convergence gradients will not be distributed as a Dirac mass any more, but the average with respect to theta will be zero (in practice small).
\end{enumerate}

\begin{figure}
	\includegraphics[width=3.5in]{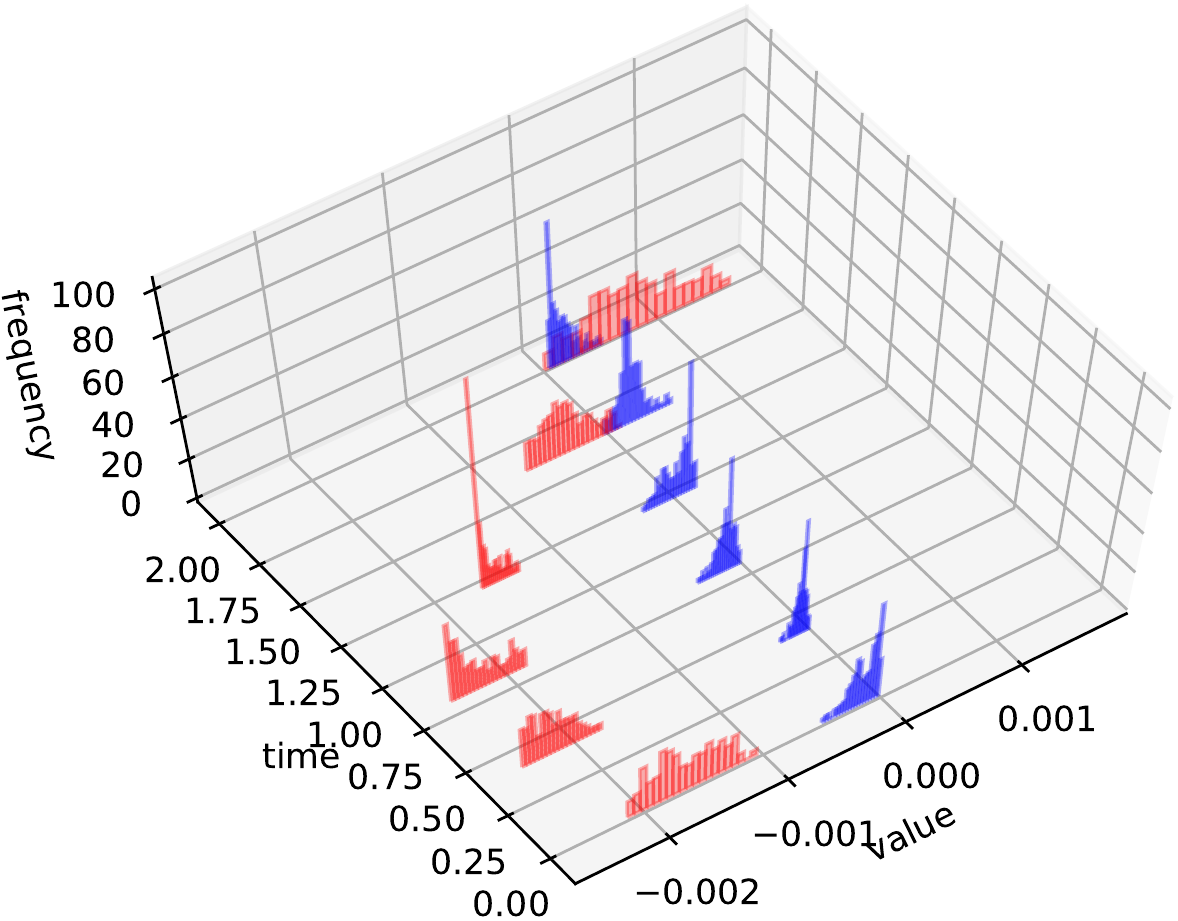}
	\caption{\label{fig:histograms1} Histogram of the gradients 
		$\nabla_{u_1(t)} {\mathcal J}(u(t_n),\theta) $
		 computed over the test sample $\theta_1^{test}$, $\theta_2^{test}$, ..., $\theta_{M_{test}}^{test}$
		 (recall $M_{test}=300$). Six time instants $t$ are
		chosen uniformly in $[0,T]$:	$t=0$, $T/5$, $2T/5$, ..., $T$.
		In red are the gradients at $u=u^1$ (iteration $k=1$) and in blue the gradients at 
		$u=u^{500}$ (iteration $k=500$).
		Here we consider the case in Section~\ref{sec:Rabitz2level}, see Figure~\ref{fig:histograms2} for the
		test case in Section~\ref{sec:Wang}.  
	}
\end{figure}

\begin{figure}
	\includegraphics[width=3.5in]{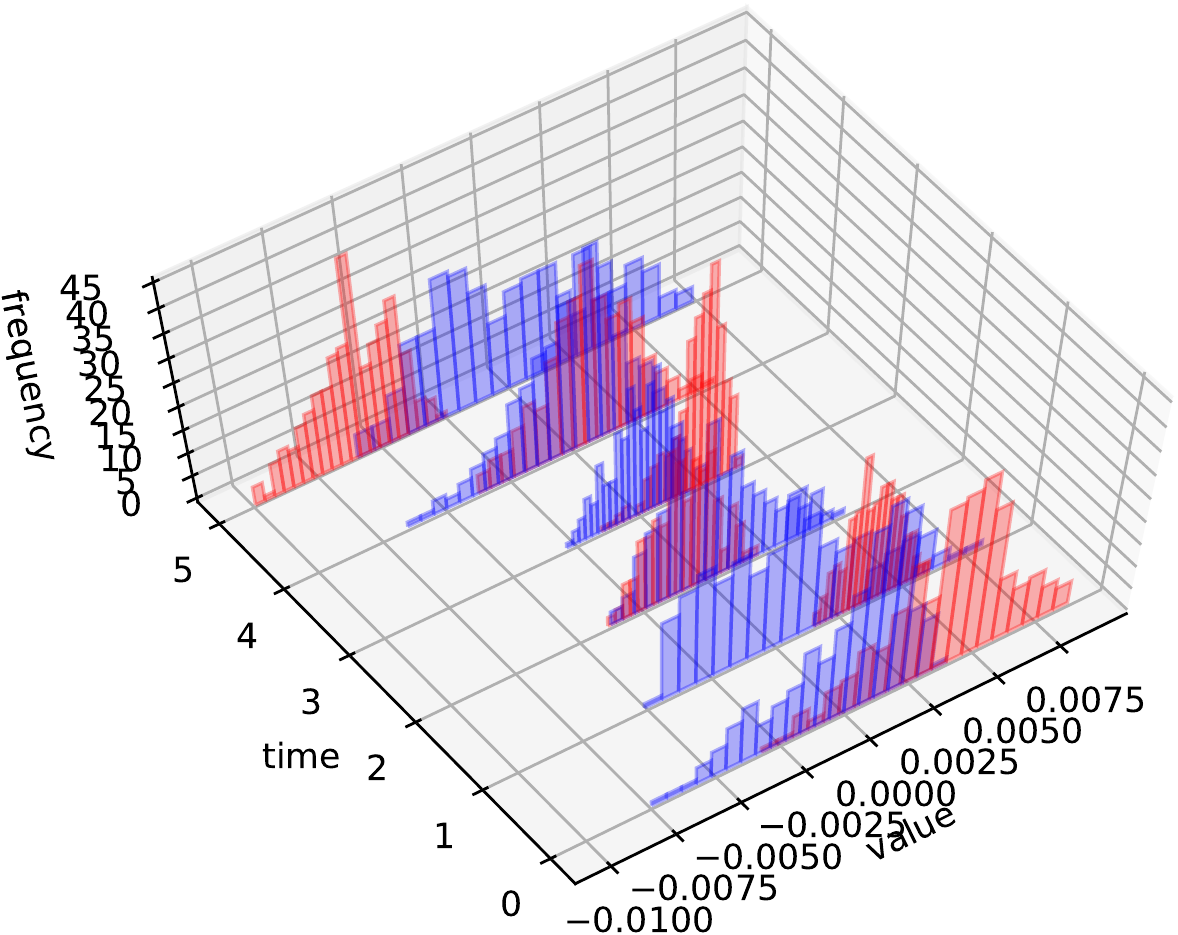}%
	\caption{\label{fig:histograms2} Histogram of the gradients 
		as in Figure~\ref{fig:histograms1} except that here the results correspond to the 
		test case in Section~\ref{sec:Wang}. 
	}
\end{figure}

We illustrate this behavior in figures \ref{fig:histograms1} and~\ref{fig:histograms2} where we plot the histograms of the gradient 
(with respect to the 
first field) 
$\nabla_{u_1(t)} {\mathcal J}(u(t_n),\theta) $ as random variables of $\theta$ 
at some time snapshots $t$. It is noticed that while in 
the first example it is possible to reduce significantly the gradient absolute value for all members of the sample (because simultaneous 
controllability holds true), in the second test case this reduction reaches a limit and the
algorithm tries instead to center the gradients on zero so that the average be as low as possible.

\section{Discussion and conclusion}

We proposed and tested in this work a stochastic approach to compute the optimal controls of inhomogeneous quantum ensembles. The algorithms have been 
employed before in other areas of stochastic optimization but not tested in this context 
(see~\cite{Rebing_robust19} for similar algorithms). Their specificity is to draw at each iteration a new set of parameters from the inhomogeneous distribution. Although at first the intuition may not recommend such an approach,
the numerical results indicate not only convergence but also faster convergence than methods based on fixed samples. In addition the method can address situations when the space of parameters is large and was tested successfully on a 6-dimensional example.

For lower dimensional examples (as in Sections ~\ref{sec:Rabitz2level} and
	~\ref{sec:Rabitz3level}) the acceleration of the stochastic
algorithms (SGD, Adam) is due essentially to the 
lower effort per iteration compared to 
a fixed grid sampling (both being proportional to the number of samples used). In higher dimensions the fixed grid approach is inherently 
less efficient due to the {\it curse of dimensionality} and may even be prohibitively large.

On the other hand, compared with SGD, the Adam algorithm has the advantage to be more robust with respect to the 
choice of the learning rate $\alpha$, but the controls are less regular.

Finally, one of the limitations of this work is to use constant learning rates. 
	Variable learning rates are potentially interesting as it could speed up convergence in the initial phases by using large values of $\alpha$ and avoid oscillations in the end by lowering $\alpha$. 
	Several schedules are proposed in the stochastic optimization literature (inverse linear, piecewise constant, ...) but their analysis remains for future work.

\appendix*
\section{Gradient computation} \label{sec:gradientcomput}

We detail below the computation of the gradient for a single  parameter $\theta$, the general case being just a mean over $\theta$.
Consider the so-called adjoint state $\lambda(t;\theta)$; it is defined at the final time as the derivative of the outcome with respect
to $C(T;\theta)$. For instance, for sections \ref{sec:Rabitz2level} - \ref{sec:Rabitz3level}: 
$\lambda(T;\theta)= -\langle C_{\text{target}}, C(T,\theta) \rangle C_{\text{target}} $ 
while for sections \ref{sec:Li} - \ref{sec:Wang} we set  $\lambda(T;\theta)=-1$.
Then for $t<T$, $\lambda(t;\theta)$ is the solution of the (backward) equation
$\frac{d}{dt} \lambda(t;\theta) = X(t,\theta)^\dagger \lambda(t;\theta)$, where 
$ X(t,\theta)^\dagger$ is the transpose conjugate of $X$ when $X$ has complex entries (examples \ref{sec:Rabitz2level} and \ref{sec:Rabitz3level})
and reduces to the transpose when $X$ is a real matrix (examples  \ref{sec:Li} and \ref{sec:Wang}).
Then  
$\nabla_{u(t)} {\mathcal J} = \langle \lambda(t;\theta), \frac{\partial X(t;\theta)}{\partial u(t)} C(t;\theta)\rangle$. In practice, given that $u$ is discretized, the state $C$ and the adjoint state $\lambda$ are also discretized at time instants
$t_n = n \Delta t$: 
$C_n(\theta) \simeq C(t_n;\theta)$,
$\lambda_n(\theta) \simeq \lambda(t_n;\theta)$ which satisfy
 $C_{n+1}(\theta)  = e^{\Delta t X(u(t_n);\theta)}C_{n}(\theta) $ and
$\lambda_n(\theta)  = e^{\Delta t X(u(t_n);\theta)^\dagger }\lambda_{n+1}(\theta) $
and the exact discrete gradient is 
$\nabla_{u(t_n)} {\mathcal J} = \langle \lambda_{n+1}(\theta), \frac{\partial  e^{\Delta t X(u(t_n);\theta)}}{\partial u(t_n)} C_n(\theta)\rangle$.

Finally, in order to compute $\frac{\partial  e^{\Delta t X(u(t_n);\theta)}}{\partial u(t_n)}$ we use  a "divide and conquer" approach coupled with a $8$-th order expansion 
as in~\cite[formula (11)]{expm8thorder17}) to obtain at the same time the exponential and the gradient (\cite[Chapter VI]{automatic_diff_rall1981}) from the knowledge of the inputs $X(u(t_n);\theta)$ and  $\frac{\partial  {X(u(t_n);\theta)}}{\partial u_k(t_n)}$. 


\begin{thebibliography}{10}
	
	\bibitem{ALLARD19987}
	Peter Allard, Magnus Helgstrand, and Torleif Hard.
	\newblock The complete homogeneous master equation for a heteronuclear two-spin
	system in the basis of cartesian product operators.
	\newblock {\em Journal of Magnetic Resonance}, 134(1):7 -- 16, 1998.
	
	\bibitem{augier}
	Nicolas {Augier}, Ugo {Boscain}, and Mario {Sigalotti}.
	\newblock {Adiabatic ensemble control of a continuum of quantum systems.}
	\newblock {\em {SIAM J. Control Optim.}}, 56(6):4045--4068, 2018.
	
	\bibitem{expm8thorder17}
	Philipp {Bader}, Sergio {Blanes}, and Fernando {Casas}.
	\newblock {An improved algorithm to compute the exponential of a matrix}.
	\newblock {\em arXiv e-prints}, page arXiv:1710.10989, Oct 2017.
	
	\bibitem{beauchard-bloch}
	Karine Beauchard, Jean-Michel Coron, and Pierre Rouchon.
	\newblock Controllability issues for continuous-spectrum systems and ensemble
	controllability of {B}loch equations.
	\newblock {\em Comm. Math. Phys.}, 296(2):525--557, 2010.
	
	\bibitem{Belhadj201523}
	Mohamed Belhadj, Julien Salomon, and Gabriel Turinici.
	\newblock Ensemble controllability and discrimination of perturbed bilinear
	control systems on connected, simple, compact {Lie} groups.
	\newblock {\em European Journal of Control}, 22(0):23 -- 29, 2015.
	
	\bibitem{borzibook}
	Alfio {Borz\`{\i}}, Gabriele {Ciaramella}, and Martin {Sprengel}.
	\newblock {\em {Formulation and numerical solution of quantum control
			problems.}}, volume~16.
	\newblock Philadelphia, PA: Society for Industrial and Applied Mathematics
	(SIAM), 2017.
	
	\bibitem{brifrabitz10}
	Constantin Brif, Raj Chakrabarti, and Herschel Rabitz.
	\newblock Control of quantum phenomena: past, present and future.
	\newblock {\em New Journal of Physics}, 12(7):075008, 2010.
	
	\bibitem{praRabitz14ensemble}
	Chunlin Chen, Daoyi Dong, Ruixing Long, Ian~R. Petersen, and Herschel~A.
	Rabitz.
	\newblock Sampling-based learning control of inhomogeneous quantum ensembles.
	\newblock {\em Phys. Rev. A}, 89:023402, Feb 2014.
	
	\bibitem{Grigoriu09}
	Jean-Michel Coron, Andreea Grigoriu, C{\u{a}}t{\u{a}}lin Lefter, and Gabriel
	Turinici.
	\newblock Quantum control design by {Lyapunov} trajectory tracking for dipole
	and polarizability coupling.
	\newblock {\em New Journal of Physics}, 11(10):105034, oct 2009.
	
	\bibitem{PhysRevA.100.022302}
	Hai-Jin Ding and Re-Bing Wu.
	\newblock Robust quantum control against clock noises in multiqubit systems.
	\newblock {\em Phys. Rev. A}, 100:022302, Aug 2019.
	
	\bibitem{Dion99}
	Claude~M. Dion, Arne Keller, Osman Atabek, and Andr\'e~D. Bandrauk.
	\newblock Laser-induced alignment dynamics of {HCN}: Roles of the permanent
	dipole moment and the polarizability.
	\newblock {\em Phys. Rev. A}, 59:1382--1391, Feb 1999.
	
	\bibitem{Glaser98}
	S.~J. Glaser, T.~Schulte-Herbr{\"u}ggen, M.~Sieveking, O.~Schedletzky, N.~C.
	Nielsen, O.~W. S{\o}rensen, and C.~Griesinger.
	\newblock Unitary control in quantum ensembles: Maximizing signal intensity in
	coherent spectroscopy.
	\newblock {\em Science}, 280(5362):421--424, 1998.
	
	\bibitem{PhysRevApplied.4.024012}
	I.~N. Hincks, C.~E. Granade, T.~W. Borneman, and D.~G. Cory.
	\newblock Controlling quantum devices with nonlinear hardware.
	\newblock {\em Phys. Rev. Applied}, 4:024012, Aug 2015.
	
	\bibitem{KhanejaPnas03}
	Navin Khaneja, Burkhard Luy, and Steffen~J. Glaser.
	\newblock Boundary of quantum evolution under decoherence.
	\newblock {\em Proceedings of the National Academy of Sciences},
	100(23):13162--13166, 2003.
	
	\bibitem{PhysRevLett.102.080501}
	Kaveh Khodjasteh and Lorenza Viola.
	\newblock Dynamically error-corrected gates for universal quantum computation.
	\newblock {\em Phys. Rev. Lett.}, 102:080501, Feb 2009.
	
	\bibitem{kingma2014adam}
	Diederik~P Kingma and Jimmy Ba.
	\newblock Adam: A method for stochastic optimization.
	\newblock {\em arXiv preprint arXiv:1412.6980}, 2014.
	\newblock {ICLR} Proceedings 2015.
	
	\bibitem{kosut}
	Robert~L. Kosut, Matthew~D. Grace, and Constantin Brif.
	\newblock Robust control of quantum gates via sequential convex programming.
	\newblock {\em Phys. Rev. A}, 88:052326, Nov 2013.
	
	\bibitem{Kuang18}
	Sen Kuang, Peng Qi, and Shuang Cong.
	\newblock Approximate time-optimal control of quantum ensembles based on
	sampling and learning.
	\newblock {\em Physics Letters A}, 382(28):1858 -- 1863, 2018.
	
	\bibitem{li-khaneja-06bis}
	J.-S. Li and N.~Khaneja.
	\newblock Control of inhomogeneous quantum ensembles.
	\newblock {\em Phys. Rev. A}, 73:030302, 2006.
	
	\bibitem{LiRuthsJCP09}
	Jr-Shin Li, Justin Ruths, and Dionisis Stefanatos.
	\newblock A pseudospectral method for optimal control of open quantum systems.
	\newblock {\em The Journal of Chemical Physics}, 131(16):164110, 2009.
	
	\bibitem{LiPNAS11}
	Jr-Shin Li, Justin Ruths, Tsyr-Yan Yu, Haribabu Arthanari, and Gerhard Wagner.
	\newblock Optimal pulse design in quantum control: A unified computational
	method.
	\newblock {\em Proceedings of the National Academy of Sciences},
	108(5):1879--1884, 2011.
	
	\bibitem{bookquantumcomputation}
	Michael~A. Nielsen and Isaac~L. Chuang.
	\newblock {\em Quantum Computation and Quantum Information: 10th Anniversary
		Edition}.
	\newblock Cambridge University Press, New York, NY, USA, 10th edition, 2011.
	
	\bibitem{nr3}
	W.H. Press, S.A. Teukolsky, W.T. Vetterling, and B.P. Flannery.
	\newblock {\em "Numerical Recipes 3rd Edition: The Art of Scientific
		Computing"}.
	\newblock Cambridge University Press, 2007.
	
	\bibitem{rabitz:hal-00536535}
	Herschel Rabitz and Gabriel Turinici.
	\newblock {Controlling quantum dynamics regardless of laser beam spatial
		profile and molecular orientation}.
	\newblock {\em Physical review A: Atomic, Molecular and Optical Physics},
	75(4):043409, 2007.
	
	\bibitem{automatic_diff_rall1981}
	Louis~B. Rall.
	\newblock {\em {Automatic Differentiation: Techniques and Applications}},
	volume 120 of {\em Lecture Notes in Computer Science}.
	\newblock Springer, Berlin, 1981.
	
	\bibitem{SGD}
	Herbert Robbins and Sutton Monro.
	\newblock A stochastic approximation method.
	\newblock {\em Ann. Math. Statistics}, 22:400--407, 1951.
	
	\bibitem{ruder2016overview}
	Sebastian Ruder.
	\newblock An overview of gradient descent optimization algorithms.
	\newblock {\em arXiv preprint arXiv:1609.04747}, 2016.
	
	\bibitem{LiJCP11}
	Justin Ruths and Jr-Shin Li.
	\newblock A multidimensional pseudospectral method for optimal control of
	quantum ensembles.
	\newblock {\em The Journal of Chemical Physics}, 134(4):044128, 2011.
	
	\bibitem{SKINNER20038}
	Thomas~E Skinner, Timo~O Reiss, Burkhard Luy, Navin Khaneja, and Steffen~J
	Glaser.
	\newblock Application of optimal control theory to the design of broadband
	excitation pulses for high-resolution nmr.
	\newblock {\em Journal of Magnetic Resonance}, 163(1):8 -- 15, 2003.
	
	\bibitem{Turinici20011}
	Gabriel Turinici and Herschel Rabitz.
	\newblock Quantum wavefunction controllability.
	\newblock {\em Chemical Physics}, 267(1--3):1 -- 9, 2001.
	
	\bibitem{0305-4470-37-1-019}
	Gabriel Turinici, Viswanath Ramakhrishna, Baiqing Li, and Herschel Rabitz.
	\newblock Optimal discrimination of multiple quantum systems: controllability
	analysis.
	\newblock {\em Journal of Physics A: Mathematical and General}, 37(1):273,
	2004.
	
	\bibitem{WangAutomaticaEnsemble18}
	Shuo Wang and Jr-Shin Li.
	\newblock Free-endpoint optimal control of inhomogeneous bilinear ensemble
	systems.
	\newblock {\em Automatica}, 95:306 -- 315, 2018.
	
	\bibitem{robust3D}
	C.~{Wu}, B.~{Qi}, C.~{Chen}, and D.~{Dong}.
	\newblock Robust learning control design for quantum unitary transformations.
	\newblock {\em IEEE Transactions on Cybernetics}, 47(12):4405--4417, Dec 2017.
	
	\bibitem{Rebing_robust19}
	Re-Bing Wu, Haijin Ding, Daoyi Dong, and Xiaoting Wang.
	\newblock Learning robust and high-precision quantum controls.
	\newblock {\em Phys. Rev. A}, 99:042327, Apr 2019.
	
	\bibitem{PhysRevA.29.1419}
	Bernard Yurke and John~S. Denker.
	\newblock Quantum network theory.
	\newblock {\em Phys. Rev. A}, 29:1419--1437, Mar 1984.
	
	\bibitem{sugny_timeoptimal}
	Y.~Zhang, M.~Lapert, D.~Sugny, M.~Braun, and S.~J. Glaser.
	\newblock Time-optimal control of spin 1/2 particles in the presence of
	radiation damping and relaxation.
	\newblock {\em The Journal of Chemical Physics}, 134(5):054103, 2011.
	
\end{thebibliography}
\bibliographystyle{plain}

\end{document}